\documentstyle{article}
\begin{document}
\begin{center}
{{\LARGE  Performance Evaluation of Switched Discrete Event
Systems\footnote[1]
 { Supported by National Natural Science Foundation of China(69925307). Email: longwang@mech.pku.edu.cn}
}}
\end{center}

\vskip 0.6cm \centerline{Long Wang} \vskip 6pt
\centerline{\small{Center for Systems and Control, Department of
Mechanics and Engineering Science }} \centerline{\small{Peking
University, Beijing 100871, CHINA }}

\vskip 0.6cm

\vskip 6pt
\begin{minipage}[t]{14cm}{
\noindent{Abstract:} This paper discusses the asymptotic periodic
behavior of a class of switched discrete event systems, and shows
how to evaluate the asymptotic performance of such systems.


}
\end{minipage}

\section{Introduction}
\par \indent

Based on max-plus algebra, a class of discrete event processes can
be described by linear recursive equations\cite{cohen}. Such a
system exhibits asymptotic periodic behavior, and its performance
can be evaluated by calculating the eigenvalue of the system
matrix in max-plus algebra. On the other hand, control techniques
based on switching among different controllers have been explored
extensively in recent years, where they have been shown to achieve
better dynamic performance\cite{morse}.

This paper proposes a new model for a class of switched discrete
event systems. Such a model consists of a finite set of discrete
event subsystems, and a switching law that orchestrates the
switching among them. We show that, the switched system can be
transformed into a non-switched system, and under certain
conditions, the switched system exhibits asymptotic periodic
behavior, and its performance can be evaluated by calculating the
eigenvalue of certain matrix in max-plus algebra.

\vskip -1in
\section{\bf{ Preliminaries}}
\par \indent

{ Denote

$$R_e = R \cup \{ - \infty \}$$

$$\epsilon = - \infty$$

\noindent and for any $x, y \in R_e$, define

$$x \oplus y = \max \{ x, y \}$$

$$x \otimes y = x + y$$

A matrix $A = (a_{ij}) \in R_e^{n \times n}$ is said to be
irreducible if $\forall i, j, \exists (i_1 = i, i_2, ......,
i_{k-1}, i_k = j)$, s.t. $a_{i i_2} + a_{i_2 i_3} + ...... +
a_{i_{k-1}j} > - \infty$.

For any matrices $A, B \in R_e^{n \times n}$, define

$$A \oplus B = ( a_{ij} \oplus b_{ij})$$

$$A \otimes B = ( \bigoplus_{k=1}^n ( a_{ik} \otimes b_{kj})) =: AB$$

Given any matrix $A \in R_e^{n \times n}$, the corresponding
directed graph (digraph) is a graph with $n$ nodes, and there is a
directed arc from node $j$ to node $i$ with weight $a_{ij}$ if and
only if $a_{ij} \neq - \infty$.

In a digraph, a circuit is a directed path that starts and ends at
the same node. In a circuit, the sum of the weights of all its
arcs divided by the number of arcs is called the mean weight. The
circuit with the maximal mean weight in a digraph is called the
critical circuit.

A zero vector is a vector with all its entries equal to $-
\infty$.

For an irreducible matrix $A \in R_e^{n \times n}$, if there exist
a real number $\lambda$ and a nonzero vector $h \in R_e^{n \times
1}$ such that $A h= \lambda h$. Then $\lambda$ and $h$ are called
the eigenvalue, eigenvector of $A$, respectively.

{\bf Lemma 1}\cite{cohen}\quad

For an irreducible matrix $A \in R_e^{n \times n}$, there is a
unique eigenvalue $\lambda$, and it equals the mean weight of the
critical circuit of its corresponding digraph.

{\bf Lemma 2}\cite{cohen}\quad

For an irreducible matrix $A \in R_e^{n \times n}$, there exist
positive integers $k_0$ and $d$ such that

$$A^{k+d} = \lambda^d A^k , \quad\quad\quad k \geq k_0$$

\noindent $d$ is called the period order of $A$.

{\bf Lemma 3}\cite{cohen}\quad

For an irreducible matrix $A \in R_e^{n \times n}$, suppose its
eigenvalue is $\lambda$, and its period order is $d$. Then there
exists a positive integer $k_0$ such that the solution of

$$X(k+1) = A X(k)$$

\noindent satisfies

$$X(k+d) = \lambda^d X(k) , \quad\quad\quad k \geq k_0$$

This shows that the system will exhibit periodic behavior
asymptotically. The mean period is exactly equal to the eigenvalue
of $A$. Hence, the eigenvalue of $A$ is an important performance
index of the system.

\vskip -1in
\section{\bf{ Switched Systems}}
\par \indent

For notational simplicity, we first discuss switching between two
subsystems\cite{morse}. That is, the switched system is governed
by

\begin{equation}
X(k+1) = A_i X(k)
\end{equation}

\noindent where $A_i \in R_e^{n \times n}$, and the switching law
is

$$i = \left\{ \begin{array}{ll}
                   1 & \quad \mbox{$k$ even}\\[3mm]
                   2 & \quad \mbox{$k$ odd}
              \end{array}
      \right. $$

\noindent Namely

$$X(1) = A_1 X(0)$$

$$X(2) = A_2 X(1)$$

$$X(3) = A_1 X(2)$$

$$X(4) = A_2 X(3)$$

$$......$$

\noindent That is

$$X(2) = A_2 X(1) = A_2 A_1 X(0)$$

$$X(4) = A_2 X(3) = A_2 A_1 X(2)$$

$$......$$

\noindent Let

\begin{equation}
Y(k) = X(2k)
\end{equation}

\noindent Then

$$Y(1) = A_2 A_1 Y(0)$$

$$Y(2) = A_2 A_1 Y(1)$$

$$......$$

\begin{equation}
Y(k+1) = A_2 A_1 Y(k)
\end{equation}

\noindent In this way, we tramsform a switched system into a
non-switched system. Thus, the following problem naturally arises:
Suppose $A_1$ and $A_2$ are irreducible matrices, is their product
$A_2 A_1$ still irreducible?

The answer is NO in general case. Consider the two irreducible
matrices

$$A_1 =   \left[ \begin{array}{ccc}
           \epsilon & 1 & \epsilon\\[3mm]
           \epsilon & \epsilon & 1\\[3mm]
           1 & \epsilon & \epsilon
           \end{array}
           \right] \quad\quad\quad\quad\quad\quad
           A_2 =   \left[ \begin{array}{ccc}
           \epsilon & \epsilon & 1\\[3mm]
           1 & \epsilon & \epsilon\\[3mm]
           \epsilon & 1 & \epsilon
           \end{array}
           \right] $$

\noindent Then, their product is

$$A_2 A_1 =   \left[ \begin{array}{ccc}
           2 & \epsilon & \epsilon\\[3mm]
           \epsilon & 2 & \epsilon\\[3mm]
           \epsilon & \epsilon & 2
           \end{array}
           \right]$$

\noindent Clearly, $A_2 A_1$ is reducible. However, if every main
diagonal entry of $A_1$ (or $A_2$) is not the null element
$\epsilon$, then the answer to the question above is YES.

{\bf Theorem 1}\quad

Suppose $A, B \in R_e^{n \times n}$ are irreducible matrices, with
all the main diagonal entries of $A$ (or $B$) not equal to
$\epsilon$. Then, $AB$ is irreducible, too.

Proof: Without loss of generality, suppose all the main diagonal
entries of $A$ are not equal to $\epsilon$. Then, for any $1 \leq
s, t \leq n $, $(AB)_{st} \neq \epsilon$ whenever $(B)_{st} \neq
\epsilon$. Moreover, since $B$ is irreducible, by definition, $AB$
is irreducible, too.

{\bf Theorem 2}\quad

Suppose $A_1, A_2 \in R_e^{n \times n}$ are irreducible matrices,
with all the main diagonal entries of $A_1$ (or $A_2$) not equal
to $\epsilon$. Then, there exist positive number $\lambda$,
positive integers $d$ and $k_0$, such that the switched system (1)
satisfies

$$X(k+d) = \lambda^d X(k) , \quad\quad\quad k \geq k_0$$

Proof: By the transformation (2), the switched system (1) can be
transformed into a non-switched system (3). That is

$$Y(k+1) = A_2 A_1 Y(k)$$

\noindent By Theorem 1, $A_2 A_1$ is irreducible. Hence, by Lemma
3 and by the transformation (2), we get the result.

{\bf Example 1}\quad

Consider the two irreducible matrices

$$A =   \left[ \begin{array}{ccc}
           2 & \epsilon & 3\\[3mm]
           6 & 2 & \epsilon\\[3mm]
           \epsilon & 4 & 3
           \end{array}
           \right] \quad\quad\quad\quad\quad\quad
           B =   \left[ \begin{array}{ccc}
           \epsilon & 3 & \epsilon\\[3mm]
           \epsilon & \epsilon & 2\\[3mm]
           4 & \epsilon & \epsilon
           \end{array}
           \right] $$

\noindent It is easy to see that

$$\lambda (A) = \frac{13}{3}, \quad\quad\quad \lambda (B) = 3$$

\noindent Moreover

$$AB =   \left[ \begin{array}{ccc}
           7 & 5 & \epsilon\\[3mm]
           \epsilon & 9 & 4\\[3mm]
           7 & \epsilon & 6
           \end{array}
           \right] $$

\noindent is also irreducible, and

$$\lambda (AB) = 9$$

\noindent Note that

$$\lambda (AB) > \lambda (A) + \lambda (B)$$

\noindent But this inequality is not always true in general case.

{\bf Example 2}\quad

Consider the two irreducible matrices

$$A =   \left[ \begin{array}{ccc}
           10 & 1 & \epsilon\\[3mm]
           \epsilon & 1 & 1\\[3mm]
           1 & \epsilon & 1
           \end{array}
           \right] \quad\quad\quad\quad\quad\quad
           B =   \left[ \begin{array}{ccc}
           1 & 1 & \epsilon\\[3mm]
           \epsilon & 1 & 1\\[3mm]
           1 & \epsilon & 10
           \end{array}
           \right] $$

\noindent It is easy to see that

$$\lambda (A) = 10, \quad\quad\quad \lambda (B) = 10$$

\noindent Moreover

$$AB =   \left[ \begin{array}{ccc}
           11 & 11 & 2\\[3mm]
           2 & 2 & 11\\[3mm]
           2 & 2 & 11
           \end{array}
           \right] $$

\noindent is also irreducible, and

$$\lambda (AB) = 11$$

\noindent Hence

$$\lambda (AB) < \lambda (A) + \lambda (B)$$

\vskip -1in
\section{\bf{ Some Extensions}}
\par \indent

More complicated switching laws can be accommodated for
performance evaluation. Suppose $A_i \in R_e^{n \times n}, i = 1,
2, ......,m$ are irreducible matrices, with all their main
diagonal entries not equal to $\epsilon$. This switched system is
governed by

\begin{equation}
X(k+1) = A_i X(k)
\end{equation}

\noindent with switching law

$$i = \left\{ \begin{array}{ll}
                   1 & \quad \quad \mbox{$k = 0, 1, 2, ......, k_1 mod(K)$}\\[3mm]
                   2 & \quad \quad \mbox{$k = k_1 +1, k_1 +2, ......, k_2 mod(K)$}\\[3mm]
                   3 & \quad \quad \mbox{$k = k_2 +1, k_2 +2, ......, k_3 mod(K)$}\\[3mm]
                   \vdots & \quad \quad \vdots\\[3mm]
                   m & \quad \quad \mbox{$k = k_{m-1} +1, k_{m-1} +2, ......, k_m mod(K)$}
              \end{array}
      \right. $$

\noindent where $K = k_m +1$.

In this case, the transformed system is

$$Y(k+1) = A_m^{k_m - k_{m-1}}......A_3^{k_3 - k_2} A_2^{k_2 - k_1} A_1^{k_1 +1} Y(k)$$

\noindent and

$$Y(k) = X(Kk)$$

Similar asymptotic periodic properties can be established as
follows.

{\bf Theorem 3}\quad

Suppose $A_i \in R_e^{n \times n}, i = 1, 2, ......,m$ are
irreducible matrices, with all their main diagonal entries not
equal to $\epsilon$. Then, for any positive integers $l_1, i = 1,
2, ......,m$, $A_m^{l_m}......A_3^{l_3} A_2^{l_2} A_1^{l_1}$ is
irreducible, too.

{\bf Theorem 4}\quad

Suppose $A_i \in R_e^{n \times n}, i = 1, 2, ......,m$ are
irreducible matrices, with all their main diagonal entries not
equal to $\epsilon$. Then, there exist positive number $\lambda$,
positive integers $d$ and $k_0$, such that the switched system (4)
satisfies

$$X(k+d) = \lambda^d X(k) , \quad\quad\quad k \geq k_0$$

Note that even if the matrix $A$ is irreducible, its power $A^l$
can be reducible for some integer $l$. For example, let

$$A =   \left[ \begin{array}{ccc}
           \epsilon & 1 & \epsilon\\[3mm]
           \epsilon & \epsilon & 1\\[3mm]
           1 & \epsilon & \epsilon
           \end{array}
           \right] $$

\noindent Then

$$A^3 =   \left[ \begin{array}{ccc}
           3 & \epsilon & \epsilon\\[3mm]
           \epsilon & 3 & \epsilon\\[3mm]
           \epsilon & \epsilon & 3
           \end{array}
           \right] $$

\noindent which is reducible. This is why we assume that all the
main diagonal entries are not equal to $\epsilon$.

\vskip -1in
\section{\bf{ Future Research}}
\par \indent

Two issues are under investigation.

\noindent 1. what is the necessary and sufficient condition for
the product of some matrices to be irreducible? In some cases,
even if each individual matrix is reducible, their product can
still be irreducible. For example

$$A =   \left[ \begin{array}{cc}
           \epsilon & 1\\[3mm]
           \epsilon & 1
           \end{array}
           \right] \quad\quad\quad\quad\quad\quad
           B =   \left[ \begin{array}{cc}
           \epsilon & \epsilon\\[3mm]
           1 & 1
           \end{array}
           \right] $$

\noindent This issue is important in performance evaluation of
switched discrete event systems.

\noindent 2. how is the eigenvalue of the product of some matrices
related to the the eigenvalue of each individual matrix? The
eigenvalue of the product of some matrices represents the
asymptotic mean period of the switched system, thereby plays an
important role in performance evaluation.

Yet another interesting research direction is to study the
asymptotic behavior of general 2-D discrete-event
systems\cite{roesser, kurek}

$$X(m+1, n+1) = A_1 X(m+1, n) \oplus A_2 X(m, n+1) \oplus A_3 X(m, n)$$

\noindent with the boundary condition

$$X(m, 0) = X_{m0}, \quad  \quad  \quad X(0, n) = X_{0n},  \quad  \quad  \quad  \quad  \quad m ,n = 0,1,2, ......$$

\noindent Under what conditions does the system exhibit periodic
behavior (with respect to $m, n$) asymptotically? and how to
evaluate its asymptotic performance?

A popular model for 2-D systems is the so-called Roesser
model\cite{roesser}

$$\left[ \begin{array}{c}
           X^h (i+1, j)\\[3mm]
           X^v (i, j+1)
           \end{array}
           \right] = \left[ \begin{array}{cc}
           A_{11} & A_{12}\\[3mm]
           A_{21} & A_{22}
           \end{array}
           \right] \left[ \begin{array}{c}
           X^h (i, j)\\[3mm]
           X^v (i, j)
           \end{array}
           \right]$$

\noindent with the boundary condition

$$X^h (0, j) = X^h_j, \quad  \quad  \quad X^v (i, 0) = X^v_i, \quad  \quad  \quad  \quad  \quad i, j = 0,1,2, ......$$

\noindent How the system (in the max-plus algebra sense) evolves
asymptotically, and how to evaluate its asymptotic performance are
the subjects of current research.

\end{document}